\documentclass[12pt]{article}
\usepackage{graphics}
\usepackage{graphicx}
\usepackage{amsmath,amsxtra,amssymb,latexsym, amscd}
\usepackage[left=3cm,right=2.5cm,top=2.5cm,bottom=2.5cm]{geometry}
\newtheorem{thm}{Theorem}

\newtheorem{prop}{Proposition}

\begin{document}
\baselineskip=17pt
\title{\bf Ruled minimal surfaces in $\Bbb R^3$ with density $e^z$\\
 \small{}}
\author{\bf  Nguyen Minh Hoang - Doan The Hieu\\
\sl  Hue Geometry Group,\\ College of Education, Hue University\\
34 Le Loi, Hue, Vietnam\\
\sl  dthehieu@yahoo.com}
 \maketitle

\begin{abstract}
 We classify ruled minimal surfaces in $\Bbb R^3$ with density $e^z.$ It is showed that there is
 no noncylindrical ruled minimal surface and there is a family of cylindrical ruled
 minimal surfaces  in $\Bbb R^3$ with density $e^z.$ It is also
 proved that all translation minimal surfaces are ruled.
 \end{abstract}
 \vskip 1cm
\section{Introduction}
Manifolds with density, a new category in geometry, appeared in many
ways in mathematics, such as quotients of Riemannian manifolds or as
Gauss space and, and can be viewed as smooth case of Gromov's
$mm$-spaces. A density on a Riemannian manifold $M^n$  is a positive
function $e^{\varphi(x)}$ used to weight volume and hypersurface
area.  Gauss space $G^n$ is a Euclidean space with Gaussian
probability density $(2\pi)^{-\frac n2}e^{-\frac{r^2}2}$ that is
very interesting to probabilists. For more details about manifolds
with density and some first results on Morgan's grand goal to
``generalize  all of Riemannian geometry to manifolds with density"
we refer the reader to \cite{mo1},  \cite{mo2}, \cite {RCBM},
\cite{Co2}. We  refer especially the reader to Chapter 18 of
Morgan's book Geometric Measure Theory (\cite{mo2}), in which the
author described general manifolds with density and their
relationship to Perelman's proof of the Poincar\'{e} conjecture.
Following Gromov (\cite[p. 213]{gr1}) the natural generalization of
the mean curvature of hypersufaces on a manifold with density
$e^{\varphi}$ is defined as
\begin{equation}      H_{\varphi}=H-\frac 1{n-1}\frac{d\varphi}{d{\bf n}}\end{equation}
and therefore, the mean curvature of a surface in $\Bbb R^3$ with
density $e^{\varphi}$ is
\begin{equation}     H_{\varphi}=H-\frac 12\cdot\frac{d\varphi}{d{\bf n}},\end{equation}
where $H$ is the Euclidean mean curvature and ${\bf n}$ is the
normal vector field of the surface. We call  $H_{\varphi}$ the mean
curvature with density or mean $\varphi$-curvature of the surface.

The literature of minimal surfaces began with Lagrange's work in
1760. Lagrange established a PDE, named after him, for the graph of
a $C^2$-function of two variables to be minimal. At that time, the
only known solution of Lagrange's equation was planes. In 1776,
Meusnier solved that equation with an additional assumption that the
level curves were straight lines and obtained the solution of ruled
minimal surface Helicoid. It is well known that, beside the trivial
case of planes, helicoid is the unique (noncylindrical) ruled
minimal surface (see \cite{babo}, \cite{ca}). In 1835, Scherk solved
Lagrange's equation for translation functions, i.e. functions of the
type $f(x,y)=g(x)+h(y)$ and discovered Scherk's minimal surfaces.

Often, locally a regular surface can be considered as the graph of
the function
$$X: U \longrightarrow \Bbb R,$$
where $U$ is a domain in $\Bbb R^2,$ in product space $\Bbb
R^2\times\Bbb R.$ In this paper, we consider ruled minimal surfaces
in space with log-linear density $\Bbb R^3=\Bbb R^2\times\Bbb
R_\varphi,$ where $\Bbb R_\varphi$ is the real line with log-linear
density $e^{\varphi}.$ The space with log-linear density $e^\varphi$
is the first non-trivial case of the new category ``manifolds with
density'' and is equivalent with the case of space with density
$e^z$ as it is showed  in section 2.  All ruled minimal surfaces in
$\Bbb R^3$ with density $e^z$ are classified. Quite different as the
classical case, it is showed that there is no noncylindrical ruled
minimal surfaces and there is a family of cylindrical ruled
 minimal surfaces. We also consider the case of translation minimal
surfaces and it is proved that, all translation minimal surfaces are
rule.

All functions in this paper are assumed belong to the class $C^2.$

\section{Minimal surfaces in spaces with densities}
From the formula of the $\varphi$-curvature, it is clear that if we
understand the geometric meaning of $\frac{d\varphi}{d{\bf n}}$ we
can discover some simple minimal surfaces in space with density. For
example, in Gauss space $G^3, \ \frac{d\varphi}{d{\bf n}}$ is the
distance from the origin to the tangent hyperplane at the
corresponding point of the surface. So it is easy to see that in
Gauss space $G^3$ (see also \cite{Co2})

\begin{enumerate}

\item  planes  have constant mean curvarture and planes
passing through the origin are minimal;
\item spheres about the origin have constant mean curvature and the
one with radius $\frac 1{\sqrt 2} is minimal;$
\item circular cylinders with the axis passing through the origin have
 constant curvature and the one with radius 1 is minimal.
\end{enumerate}

Let $\varphi(x_1,x_2,\ldots, x_n) $ be a linear function on
Euclidean space $\Bbb R^n$
$$\varphi(x_1,x_2,\ldots, x_n)=\sum_{i=1}^n a_ix_i.$$
  A linear density on $\Bbb R^n$ is the positive function $e^{\varphi(x)},$ where
  $x=(x_1,x_2,\ldots, x_n).$ It is easy to see that, the sets of points in
  $\Bbb R^n$ with the linear density $e^\varphi$ that have the same density are hyperplanes.
   By a suitable changing the coodinate system we can assume that the density has the
   form $e^{x_n}$ and therefore we can view the space $\Bbb R^n$ with density
    $e^{\varphi(x)}$  as product $\Bbb R^{n-1}\oplus\Bbb R_\varphi,$ where $\Bbb R^{n-1}$
    is nothing but  Euclidean $(n-1)$-space and $\Bbb R_\varphi$  is a real line with
     density $e^{x_n}.$

Since $\nabla \varphi=(0,0,\ldots, 1),\ \  \frac{d\varphi}{d{\bf n}}=\langle \nabla\varphi,
 {\bf n}\rangle$ is the cosine of the angle between {\bf n} and $z$-axis.
 By the definition of the mean $\varphi$-curvature,
 it is easy to see that $H_\varphi$ does not change under
 a translation or a rotation about $z$-axis and moreover, we have

\begin{enumerate}
\item  hyperplanes in $\Bbb R^n$ with density $e^{x_n}$  have
constant mean curvarture.
\item hyperplanes in $\Bbb R^n$ with density $e^{x_n}$  that are parallel
to the $x_n$-axis  have zero mean  curvarture.
\item a circular hypercylinder with rulings parallel to $x_n$-axis have constant mean
 curvature.
   \end{enumerate}
\section{Ruled minimal surfaces in $\Bbb R^3$ with linear density $e^z$}

Now we consider the problem of classifying all ruled minimal
surfaces in $\Bbb R^3$ with a linear density. Coordinates in $\Bbb
R^3$ are denoted by $(x,y,z).$ Without loss of generality we can
assume that the density is $e^z.$

Locally, a  ruled surface is given by the equation
\begin{equation}  \label{1}    X(u,v)=\alpha(u)+v\beta(u),\ \ u\in(a,b),\ \ v\in(c,d).\end{equation}
We can assume that $|\alpha'|=1,\ \ |\beta|=1,\ \ \langle\alpha',
\beta\rangle=0.$

We will focus on two cases: cylindrical ruled surface
($\beta=\text{const.},$ for all $v\in (c,d))$ and noncylindrical
ruled surfaces ( $\beta'\ne 0$ for all $v\in (c,d)).$ Our results
show that it is no need to consider at the isolated points where
$\beta'(v)=0.$

Denote $E,F,G$ the coefficients of the first fundamental form and $e,f,g$
 the coefficients of the second fundamental form, a direct computation yields
$$ N=\frac{(\alpha'+v\beta')\wedge\beta}{|(\alpha'+v\beta')\wedge\beta)|};$$
$$E=1+2v\langle\alpha',\beta'\rangle+v^2|\beta'|^2,\ \ \ F=0,\ \ \ G=1;$$
$$e=\langle N,\alpha''+v\beta''\rangle,\ \ \ f=\langle N,\beta'\rangle,\ \ \ g=0;$$
\begin{equation}  \label {31}  H_\varphi=\frac 12\left[\frac{\langle N,\alpha''+v\beta''\rangle}{1+2v\langle\alpha',\beta'\rangle+v^2|\beta'|^2}-\langle N,\nabla\varphi\rangle\right].   \end{equation}
We have
\begin{prop} \begin{equation} \label{3}  H_\varphi=0\Leftrightarrow \begin{cases}
\langle\alpha'\wedge\beta,\alpha''\rangle=\langle\alpha'\wedge\beta, \nabla\varphi\rangle\\
\langle\alpha'\wedge\beta,\beta''\rangle+\langle\beta'\wedge\beta,\alpha''\rangle=\langle\beta'\wedge\beta, \nabla\varphi\rangle+\langle\alpha'\wedge\beta, 2\langle\alpha',\beta'\rangle\nabla\varphi\rangle\\
\langle\beta'\wedge\beta,\beta''\rangle=\langle\beta'\wedge\beta, 2\langle\alpha',\beta'\rangle\nabla\varphi\rangle+\langle\alpha'\wedge\beta, |\beta'|^2\nabla\varphi\rangle\\
\langle\beta'\wedge\beta, |\beta'|^2\nabla\varphi\rangle=0
                  \end{cases}.   \end{equation}
\end{prop}
{\bf Proof.} $H_\varphi=0$ if and only if
\begin{equation}  \label{6a}  \frac{\langle N,\alpha''+v\beta''\rangle}{1+2v\langle\alpha',\beta'
\rangle+v^2|\beta'|^2}=\langle N,\nabla\varphi\rangle   \end{equation}
Replace $ N=\frac{(\alpha'+v\beta')\wedge\beta}{|(\alpha'+v\beta')\wedge\beta)|}$ in (\ref{6a})
 we get an equality with both RHS and LHS are polynomials on the variable $v.$ Identifying
 the coefficients we obtain  (\ref{3}).\hfill$\Box$

\subsection{The case of $\beta'\ne 0$}
From the last equation of (\ref{3})
\begin{equation}  \langle\beta'\wedge\beta, |\beta'|^2\nabla\varphi\rangle=0,    \end{equation}  we have  $\langle\beta'\wedge\beta, \nabla\varphi\rangle=0.$
Note that $\nabla\varphi=(0,0,1),$ and $\beta\perp\beta',$ we conclude $\beta $
belongs to a plane containing the $z$-axis. After a rotation about the $z$-axis,
 we can assume that $\beta=(\cos t(u),0,\sin t(u)),$ with $t'\ne 0$ and therefore, the third equality
  of (\ref{3}) becomes
\begin{equation}\label{5}   \langle\alpha'\wedge\beta, |\beta'|^2\nabla\varphi\rangle=0.    \end{equation}

From (\ref{5}) we conclude that $\alpha'$ belongs to the plane $\{y=0\}$ and the curve $\alpha$ lies on a plane parallel to $xz$-plane. It is clear that, $\alpha$ and $\beta$ satisfy system of equations (\ref{3}). Thus,

\begin{prop} If $\beta'\ne 0, \ \forall v\in(c,d),$ ruled minimal surfaces determined
by (\ref{1}) are planes parallel to the $z$-axis.
\end{prop}

\subsection{The case of $\beta'=0$}

Since $\beta'=0,\ \beta =(a,b,c)=\text{const.}$ and $\ a^2+b^2+c^2=1,$  system (\ref{3})
becomes
\begin{equation}  \label{7}\begin{cases}
\langle\alpha'\wedge\beta, \alpha''\rangle=\langle\alpha'\wedge\beta, \nabla\varphi\rangle\\
\beta=(a,b,c)=\text{const.}
\end{cases}.    \end{equation}

The first equation of (\ref{7}) implies
$$\alpha''-\nabla\varphi=m\alpha'+n\beta,$$
and hence
$$\begin{aligned}   \langle\alpha''-\nabla\varphi, \alpha'\rangle&=m=-\langle\nabla\varphi,\alpha'\rangle,\\
\langle\alpha''-\nabla\varphi, \beta\rangle&=n=-\langle\nabla\varphi,\beta\rangle.
              \end{aligned}$$
Thus,
$$\alpha''-\nabla\varphi=-\langle\nabla\varphi,\alpha'\rangle\alpha'-\langle\nabla\varphi,\beta\rangle\beta,$$
or
\begin{equation}\label{8}  \alpha''+\langle\nabla\varphi,\alpha'\rangle\alpha'=\nabla\varphi-\langle\nabla\varphi,\beta\rangle\beta    \end{equation}

Since the mean $\varphi$- curvarture does not change under a
rotation about $z$-axis, we can assume $a=0.$  So we have
$$\nabla\varphi-\langle\nabla\varphi,\beta\rangle\beta =(0,-cb,1-c^2).$$
Because $b^2=1-c^2,$
\begin{equation} \label{8a}   \nabla\varphi-\langle\nabla\varphi,\beta\rangle\beta =(0,-cb,b^2).
\end{equation}
Suppose that $\alpha=(x(u),y(u),z(u)),$ then (\ref{8a}) is equivalent to the
 following system of equations
\begin{equation} \label{10} \begin{cases}
x''+x'z'=0,\\
y''+y'z'=-cb,\\
z''+z'^2=b^2.
                   \end{cases}    \end{equation}
We treat the two special case $\beta=(0,0,\pm1)$ and $\beta=(0,\pm1,0)$ first.

If $\beta=(0,0,\pm 1),$ the RHS of (\ref{8}) equals zero. We conclude that $\alpha'$ and $\alpha''$ are parallel, and hence $\alpha''=0.$ Thus, $\alpha$ is a straight line and we have
\begin{prop} If $\beta=(0,0,\pm 1),$ ruled minimal surfaces determined by
equation (\ref{1}) are planes parallel to the $z$-axis.
\end{prop}
We can also treat this case by solving (\ref{10}). In this case (\ref{10}) becomes
\begin{equation} \label{11} \begin{cases}
x''+x'z'&=0,\\
y''+y'z'&=0,\\
z''+z'^2&=0.
                   \end{cases}    \end{equation}
Since $\beta\perp\alpha',$ we get
$z'=0,$ and hence $x''=y''=0.$ We conclude that $\alpha=(x,y,z)$ is a straight line lying
on the plane $z=\text{const.}$ and hence the ruled surfaces is a plane parallel to the $z$-axis.

If $\beta=(0,\pm 1, 0),$ (\ref{10}) becomes
\begin{equation} \label{12} \begin{cases}
x''+x'z'&=0,\\
y''+y'z'&=0,\\
z''+z'^2&=1.
                   \end{cases}    \end{equation}
Since $\beta\perp\alpha',$ we get $y'=0$ and conclude that $\alpha$ lies on the plane
 $y=\text{const.}.$ The last equation of (\ref{12}) gives us the solution

$$z'=1-\frac {2}{1+Ae^{2u}}=\frac {Ae^{2u}-1}{Ae^{2u}+1},\ \ \ A>0$$
$$z=\log (1+Ae^{2u})-u=\log(e^{-u}+Ae^{u}).$$

The first equation in (\ref{12}) gives
$$x'=Be^{-z}\frac B{e^{-u}+Ae^u}=\frac {Be^u}{1+Ae^{2u}},$$
and hence
$$x=\frac B{\sqrt{A}}\arctan\sqrt{A}e^u+C.$$
Note that, since the mean $\varphi$-curvature does not change under a translation, we can take the constant in the expression of $x$ and $z$ to be zero.

Since $x'^2+y'^2+z'^2=1,$ we have $4A=B^2.$ Thus,
\begin{prop} If $\beta=(0,\pm 1, 0),$ a ruled minimal surface determined by equation (\ref{1})
 has a parametric equation of the following form
\begin{equation}\begin{cases}
x&=2\arctan\sqrt{A}e^u\\
y&=\pm v\\
z&=\log(e^{-u}+Ae^{u})
       \end{cases}     \end{equation}
\end{prop}
\begin{figure}
 \begin{center}
  \includegraphics[width=10cm]{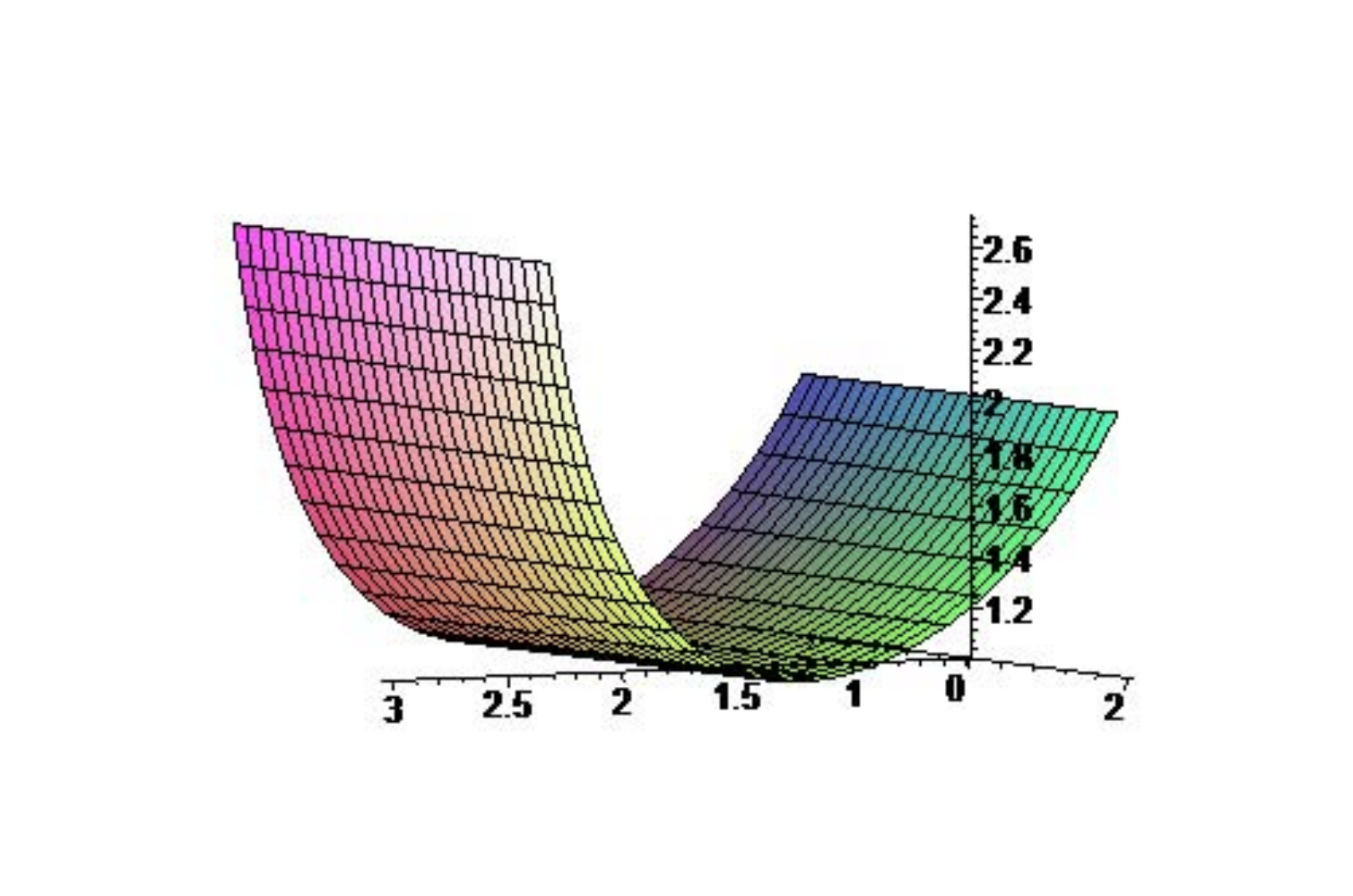}\\
  \end{center}
 \caption{Ruled minimal surface with $\beta'=(0,1,0).$}\label{h1}
   \end{figure}
\begin{figure}
 \begin{center}
  \includegraphics[width=10cm]{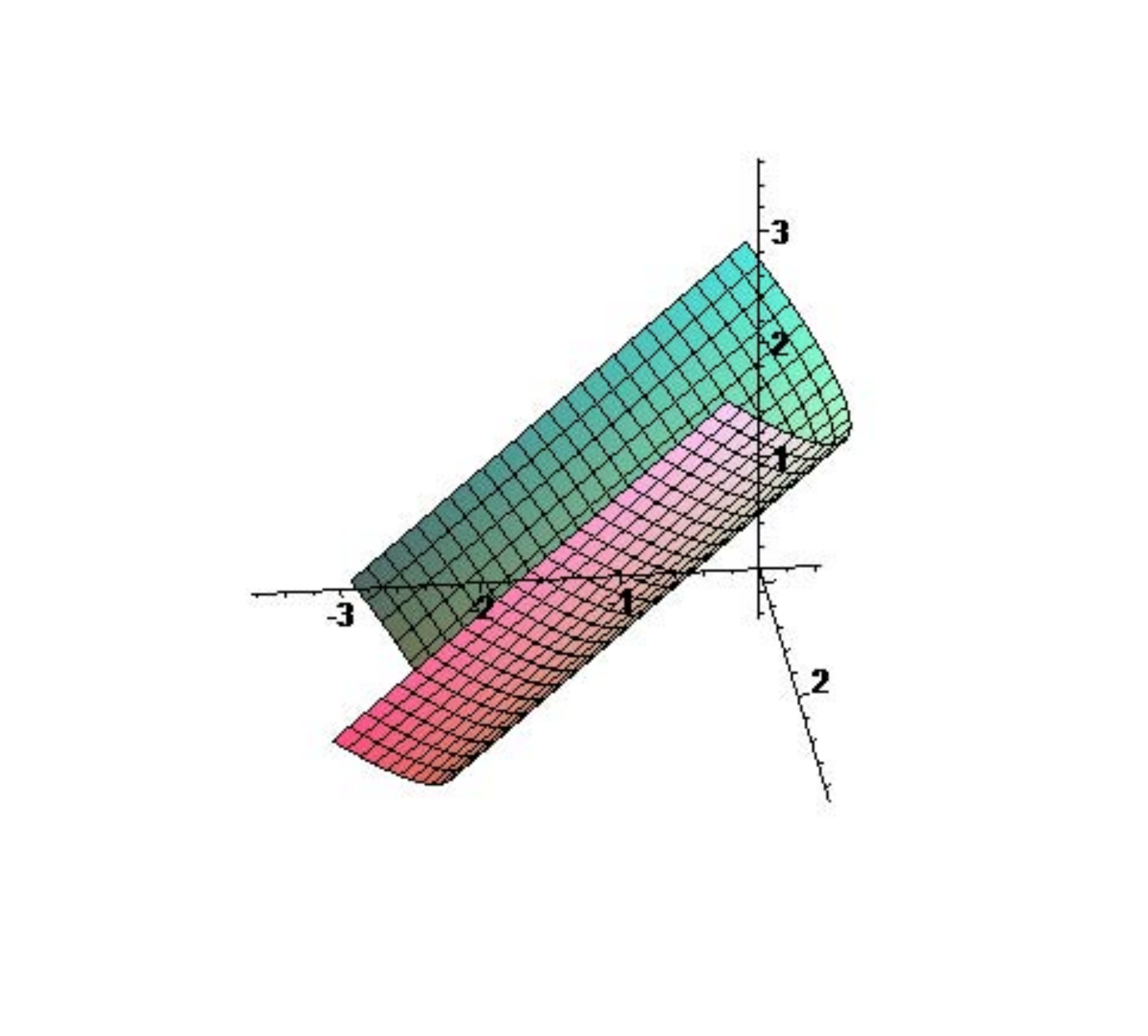}\\
  \end{center}
 \caption{Ruled minimal surface with $\beta'=(0,\frac{\sqrt2}2,\frac{\sqrt2}2).$}\label{h2}
   \end{figure}
\begin{figure}
 \begin{center}
  \includegraphics[width=10cm]{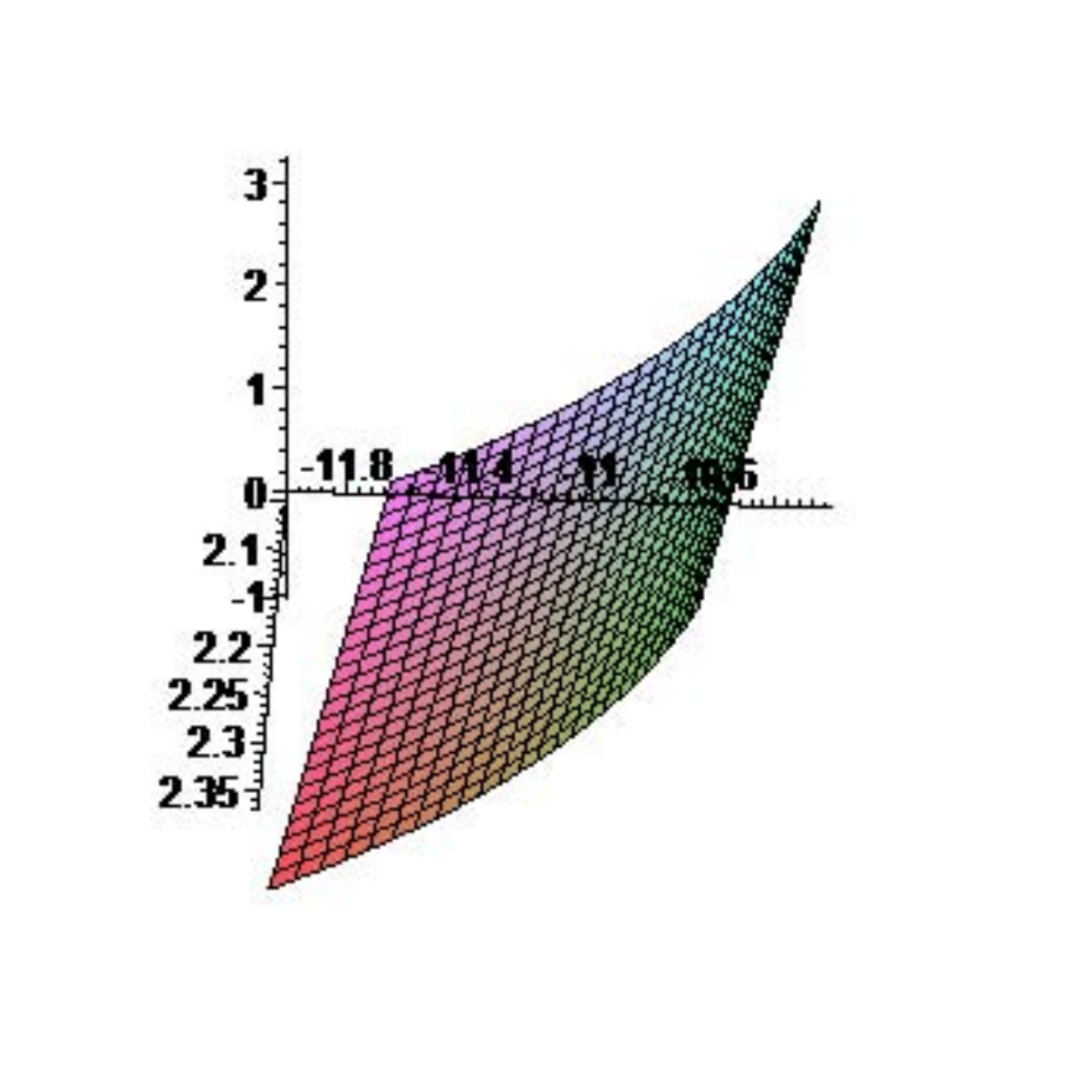}\\
  \end{center}
 \caption{Ruled minimal surface with $\beta'=(0,0.1,0.99).$}\label{h3}
   \end{figure}
\begin{figure}
 \begin{center}
  \includegraphics[width=10cm]{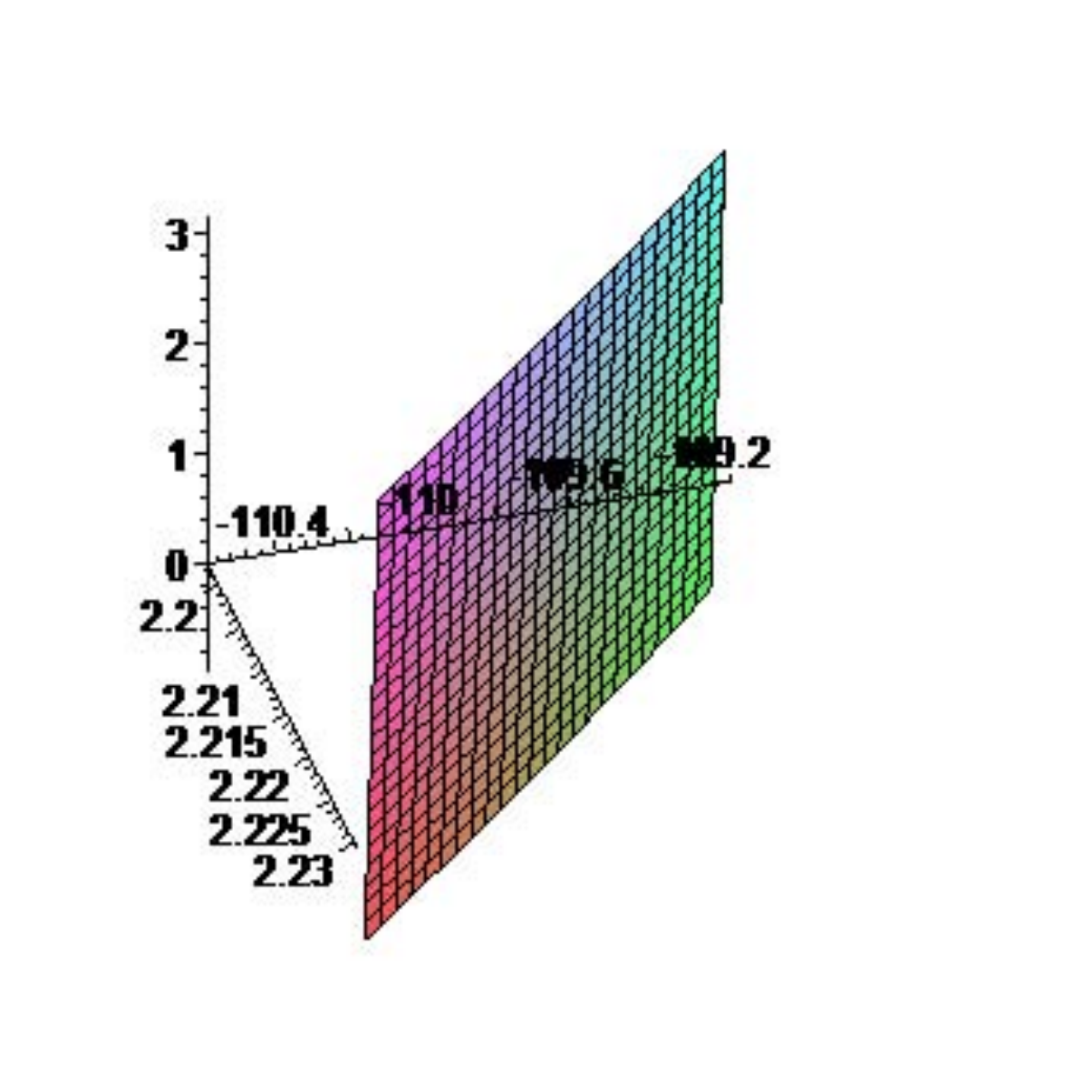}\\
  \end{center}
 \caption{Ruled minimal surface with $\beta'=(0,0.01,0.999).$}\label{h4}
   \end{figure}


If $\beta=(0,b, c),\ b,c\ne 0$ the system (\ref{10}) becomes
\begin{equation} \label{14} \begin{cases}
x''+x'z'&=0,\\
y''+y'z'&=-cb,\\
z''+z'^2&=b^2.
                   \end{cases}    \end{equation}
Since $\beta\perp\alpha',$ we get $by'=-cz'$ and conclude that $\alpha$ lies on the plane
 $by+cz+d=0.$
The last equation of (\ref{14}) gives us the solution

$$z'=b-\frac {2b}{1+Ae^{2bu}}=\frac {bAe^{2bu}-b}{Ae^{2bu}+1},\ \ \ A>0$$
$$z=\log (1+Ae^{2bu})-bu=\log(e^{-bu}+Ae^{bu}).$$

The first equation in (\ref{14}) gives
$$x'=Be^{-z}=\frac B{e^{-bu}+Ae^{bu}},$$
and hence
$$x=\frac B{b\sqrt{A}}\arctan(\sqrt{A}e^{bu})+C.$$

Since $x'^2+y'^2+z'^2=1,$ we have $4A=B^2.$ Thus,
\begin{prop} If $\beta=(0,b,c),\ b,c\ne 0,$ a ruled minimal surface determined by
equation (\ref{1}) has a parametric equation of the following form
\begin{equation}\label{17}\begin{cases}
x&= 2\arctan\sqrt{A}e^{bu},\\
y&=-\frac cb\log(e^{-bu}+Ae^{bu})+bv,\\
z&=\log(e^{-bu}+Ae^{bu})+cv.
         \end{cases}   \end{equation}

\end{prop}

Combining these above results, we have
\begin{thm}Beside planes parallel to the $z$-axis the only cylindrical ruled minimal surfaces in $\Bbb R^3$ with density $e^z$ are those given by  (\ref{17}).
\end{thm}
\section{Translation minimal surfaces in space with linear density
$e^z.$} In this section we study translation minimal surfaces in
$\Bbb R^3$ with density $e^z.$ We prove that all translation
surfaces that are minimal must be ruled.
\begin{thm}A translation surface given by
$$X(u,v)=(u,v,g(u)+h(v))$$
is minimal if either $g(u)=au+b$ or $h(v)=cv+d.$
\end{thm}
{\bf Proof.} A straightforward computation shows that
\begin{equation} \label{eq1}   H_\varphi=0\Leftrightarrow
g''(1+h'^2)+h''(1+g'^2)=1+g'^2+h'^2.
 \end{equation}
We fix $v=v_0,$ and set $A=1-h''(v_0),\ \ B=1+h'^2(v_0),\ \
C=1+h'^2(v_0)-h''(v_0).$ Note that $B>0$ and $C=B-A-1.$ Thus, $f$
must be satisfies the following equation
$$Ag'^2+Bg''=C,$$
and hence
$$g''=\frac{C-Ag'^2}B.$$

Substitute $g''$ into (\ref{eq1}), we get
\begin{equation}\label{eq2}
 g'^2\left[h''-\frac AB(1+h'^2)-1\right]=1-h''+h'^2-\frac
CB(1+h'^2).\end{equation}
 From equation (\ref{eq2}), unless
$g'=\text{const.}$ we must have
\begin{equation}\label{eq3} h''-\frac AB(1+h'^2)-1=0  \end{equation}
and
 \begin{equation} \label{eq4}    1-h''+h'^2-\frac CB(1+h'^2)=0.  \end{equation}
Substitute $h''$ from equation (\ref{eq3}) into equation (\ref{eq4})
we obtain
$$(1+h'^2)\left(1-\frac CB-\frac AB\right)=1.$$
Note that $C=B-A-1,$ we get
$$h'^2=B-1.$$
Thus, $h'=\text{const.}$ and the proof is completed.\hfill $\Box$

Since the role of $g$ and $h$ are the same, we onle need to consider
 translastion minimal surfaces of the following form
 $$X(u,v)=(u,v, g(u)+cv+d).$$
A straightforward computation shows that $g$ must be of the form
$$ -(1+c^2)\log\left|\cos\frac{u+D}{\sqrt{1+c^2}}\right|.$$

\end{document}